\documentclass[12pt]{amsart}

\usepackage{environ, xcolor}
\NewEnviron{commentA}{}

\newcommand\Aoff{\RenewEnviron{commentA}{}}
\Aoff{} 

\usepackage{amsmath, amssymb}
\usepackage{array}
\usepackage[frame,cmtip,arrow,matrix,line,graph,curve]{xy}
\usepackage{graphpap, color, paralist, pstricks}
\usepackage[mathscr]{eucal}
\usepackage[pdftex]{graphicx}
\usepackage[pdftex,colorlinks,backref=page,citecolor=blue]{hyperref}
\usepackage{cleveref}
\usepackage{tikz, verbatim}

\newtheorem{theorem}{Theorem}[section]
\newtheorem{proposition}[theorem]{Proposition}
\newtheorem{corollary}[theorem]{Corollary}
\newtheorem{lemma}[theorem]{Lemma}

\theoremstyle{definition}
\newtheorem{definition}[theorem]{Definition}
\newtheorem{example}[theorem]{Example}

\newtheorem{question}[theorem]{Question}
\newtheorem{remark}[theorem]{Remark}

\usepackage[margin=1in]{geometry} 
\usepackage{amsmath,amssymb,mathtools, mathrsfs,esint,tikz-cd,verbatim}

\newcommand{\Z}{{\mathbb Z}}
\newcommand{\Q}{{\mathbb Q}}

\newcommand{\C}{{\mathbb C}}

\newcommand{\on}[1]{\operatorname{#1}}

\newcommand{\Spec}{{\on{Spec}}}

\newcommand{\Pic}{{\on{Pic}}}

\keywords{weighted hypersurfaces, automorphism groups, Jordan constants}

\subjclass[2020]{14J50, 14J70}

\title{Automorphisms of weighted projective hypersurfaces}
\author{Louis Esser}

\address{Department of Mathematics, Princeton University, Fine Hall, Washington Road, Princeton, NJ 08544-1000, USA}
\email{esserl@math.princeton.edu}

\begin{document}

\begin{abstract}
We prove several results concerning automorphism groups of quasismooth complex weighted projective hypersurfaces; these generalize and strengthen existing results for hypersurfaces in ordinary projective space.  First, we prove in most cases that automorphisms extend to the ambient weighted projective space. We next provide a characterization of when the linear automorphism group is finite and find an explicit uniform upper bound on the size of this group.  Finally, we describe the automorphisms of a generic quasismooth hypersurface with given weights and degree.
\end{abstract}
\maketitle

\section{Introduction}

Hypersurfaces in projective space are among the most well-studied varieties.  In particular, a great deal is known about their automorphism groups, due to landmark theorems of Grothendieck-Lefschetz \cite{SGA2}, Matsumura-Monsky \cite{MM}, and others.  In this paper, we generalize and strengthen several of these results to hypersurfaces in any weighted projective space over $\C$.  

Given a collection of positive integers $a_0,\ldots,a_{n+1}$, the {\it weighted projective space} $\mathbb{P} \coloneqq \mathbb{P}(a_0,\ldots,a_{n+1})$ is defined to be the projective quotient variety $(\mathbb{A}^{n+2} \setminus \{ 0 \})/\mathbb{G}_{\on{m}}$.  Here the multiplicative group $\mathbb{G}_{\on{m}}$ acts by the formula $t \cdot (x_0,\ldots,x_{n+1}) = (t^{a_0}x_0,\ldots,t^{a_{n+1}}x_{n+1})$.  When all weights $a_i$ are equal to $1$, we recover ordinary projective space $\mathbb{P}^{n+1}$.  Unlike $\mathbb{P}^{n+1}$, weighted projective spaces are usually singular.

Hypersurfaces in weighted projective space are an extremely useful class of algebraic varieties. Indeed, many of their properties are determined combinatorially by the choice of degree and weights, and they exhibit a large range of behavior.  In particular, there is significant evidence that weighted projective hypersurfaces are flexible enough to solve a diverse range of optimization problems in algebraic geometry (see, for example, \cite{ETW,CEW,ETWindex,Totaro}).  It is therefore desirable to understand basic properties of their automorphism groups.

The outline of the paper is as follows: in \Cref{background}, we discuss the necessarily preliminaries on the objects of study. In \Cref{linearity}, we prove that all automorphisms of well-formed quasismooth weighted projective hypersurfaces $X \subset \mathbb{P}$ extend to the ambient weighted projective space if $\dim(X) \geq 3$ or $\dim(X) = 2$ and $X$ has non-trivial canonical class (\Cref{main-theorem}).  This generalizes the same statement for hypersurfaces in $\mathbb{P}^{n+1}$, which is a consequence of the Grothendieck-Lefschetz theorem. Some partial results in this direction have appeared in works of Przyjalkowski and Shramov \cite{PrS,PrS2}.  An automorphism of $X$ that extends to $\mathbb{P}$ is called {\it linear}.

We prove several results on the size of the linear automorphism group $\mathrm{Lin}(X)$ of a well-formed quasismooth $X$ in \Cref{bounds}.  In particular, we give an exact characterization of when this group is finite in terms of the degree $d$ and weights $a_0,\ldots,a_{n+1}$ of $X$ (\Cref{finite}).  When it is finite, we prove that there is an upper bound
$$|\mathrm{Lin}(X)| \leq C_n \frac{d^{n+1}}{a_0 \cdots a_{n+1}},$$
where $C_n$ depends only on the dimension $n$ (\Cref{bound}). We give a procedure for effectively calculating a value $C_n$ for which the inequality above holds using the Jordan constants of general linear groups over $\C$. We also compute the optimal value for $C_1$.

Finally, in \Cref{verygeneral}, we consider automorphisms of a very general quasismooth hypersurface $X$ with a given degree and weights.  We prove that when $d \geq 5 \max\{a_0,\ldots,a_{n+1}\}$, the group $\mathrm{Lin}(X)$ is contained in the center of $\mathrm{Aut}(\mathbb{P})$ (\Cref{generic}).  When $\mathbb{P} = \mathbb{P}^{n+1}$, the center of $\mathrm{Aut}(\mathbb{P}^{n+1}) = \mathrm{PGL}_{n+2}$ is trivial, so we recover the usual statement that $\mathrm{Lin}(X) = 1$ for $X$ very general.  However, this stronger statement is not always true for other choices of weights and degree.  We give several examples illustrating new phenomena that arise for generic automorphisms of weighted projective hypersurfaces that don't occur in ordinary projective space.

\noindent{\it Acknowledgements. }The author was
supported by NSF grant DMS-2054553. Thanks to Alex Duncan, Victor Przyjalkowski, and Burt Totaro for useful conversations.

\subsection{Background on Weighted Projective Hypersurfaces} \label{background}

Throughout the paper, we'll work over the complex numbers, though some of the statements below remain true over other fields.  For a complete introduction to weighted projective hypersurfaces, see \cite{Iano-Fletcher}.

We say that $\mathbb{P}$ is {\it well-formed} if the action of $\mathbb{G}_{\on{m}}$ on $\mathbb{A}^{n+2}$ has trivial stabilizers in codimension $1$.  This holds exactly when $\gcd(a_0,\ldots,\widehat{a_i},\ldots,a_{n+1}) = 1$ for each $i$.  Every weighted projective space is isomorphic as a variety to a well-formed one, so we will only consider well-formed spaces $\mathbb{P}$.  If $S = \C[x_0,\ldots,x_{n+1}]$ is the graded polynomial ring with the weight of each $x_i$ equal to $a_i$, then $\mathrm{Proj}(S) \cong \mathbb{P}(a_0,\ldots,a_{n+1})$.  The space $\mathbb{P}$ is equipped with a reflexive sheaf $\mathcal{O}(d)$ for every integer $d \in \Z$ from the Proj construction, which is the sheaf associated to a Weil divisor of degree $d$.  It is a line bundle if and only if $d$ is divisible by every weight $a_i$.

Let $f$ be a weighted homogeneous polynomial of degree $d$.  Then $f$ defines a hypersurface $X = \{f = 0\} \subset \mathbb{P}(a_0,\ldots,a_{n+1})$ of dimension $n$.  The {\it affine cone} $C_X$ over $X$ is the subvariety of $\mathbb{A}^{n+2}$ defined by the same equation $f$, while the {\it punctured} affine cone is $C_X^* \coloneqq C_X \setminus \{0\}$.  A hypersurface $X$ is {\it quasismooth} if the punctured affine cone is smooth, and it is {\it well-formed} if $\mathbb{P}$ is well-formed and $X \cap \mathbb{P}_{\on{sing}}$ has codimension at least $2$ in $X$.   When $X$ is well-formed and quasismooth, the canonical divisor satisfies the adjunction formula $K_X = \mathcal{O}_X(d-a_0 - \cdots -a_{n+1})$.  (A slightly weaker condition is used for well-formedness of weighted complete intersections in \cite{Iano-Fletcher}, but the stronger condition used here is the one that implies the adjunction formula, see \cite{Przy}; in any case, the two definitions coincide for hypersurfaces.) 

Quasismooth hypersurfaces of degree $d$ in $\mathbb{P}$ only exist for certain choices of weights and degree, according to the following criterion, which holds in characteristic $0$:

\begin{proposition}[{{{\cite[Theorem 8.1]{Iano-Fletcher}}}}] \label{qsmooth}
There exists a quasismooth hypersurface $X$ of degree $d$ in the weighted projective space $\mathbb{P}(a_0,\ldots, a_{n+1})$ if and only if one of the following properties holds:
\begin{enumerate}
    \item $a_i = d$ for some $i \in \{0,\ldots,n+1\}$, or
    \item for each nonempty subset $I$ of $\{0,\ldots ,n+1\}$, either
    \begin{enumerate}
        \item $d$ is an $\mathbb{N}$-linear combination of the weights $a_i$ for $i \in I$, or
        \item there are at least $|I|$ indices $j \notin I$ such that $d-a_j$ is an $\mathbb{N}$-linear combination of the numbers $a_i$ with $i \in I$.
    \end{enumerate}
\end{enumerate}
If (1) or (2) holds, then the general hypersurface of degree $d$ is quasismooth.
\end{proposition}

Here, ``$\mathbb{N}$-linear combination" means a linear combination with nonnegative integer coefficients. We'll frequently use the following version of the proposition applied to singleton sets $I = \{ i \}$: 

\begin{proposition} \label{monomialexistence}
Suppose $X = \{f = 0\} \subset \mathbb{P}$ is a quasismooth hypersurface of degree $d$ in $\mathbb{P}$. Then for each $i = 0,\ldots,n+1$, there exists a monomial of degree $d$ with nonzero coefficient in $f$ having the form either (a) $x_i^k$, or (b) $x_i^k x_j$ for some $j \neq i$.
\end{proposition}

\begin{proof}
If such a monomial does not exist for some $i$, then the function $f$ and all its derivatives would vanish at the coordinate point $e_i \in \mathbb{A}^{n+2}$, contradicting the smoothness of the punctured affine cone $C_X^*$.
\end{proof}

As mentioned above, $\mathbb{P}(a_0,\ldots,a_{n+1}) = \mathrm{Proj}(S)$, where $S = \C[x_0,\ldots,x_{n+1}]$ is the graded polynomial ring with variables of the given weights.  We also use the notation $S_m$ to denote the $m$th graded piece of $S$.  Any graded automorphism of $S$ induces an automorphism of weighted projective space; let $\mathrm{Aut}(S)$ denote the group of graded automorphisms.  In fact, every automorphism of $\mathbb{P}$ comes from $\mathrm{Aut}(S)$. This is proven in \cite[Section 8]{Amrani}, but for completeness, we provide a proof here:

\begin{lemma} \label{autP}
Let $\mathbb{P} \coloneqq \mathbb{P}(a_0,\ldots,a_{n+1})$ be a well-formed weighted projective space and $S \coloneqq \C[x_0,\ldots,x_{n+1}]$ the graded polynomial ring with the weight of $x_i$ equal to $a_i$.  Then the natural map $\mathrm{Aut}(S) \rightarrow \mathrm{Aut}(\mathbb{P})$ is surjective.  The kernel is isomorphic to $\C^*$, where the isomorphism associates to any $t \in \C^*$ the automorphism mapping $x_i \mapsto t^{a_i} x_i$ for each $i$.
\end{lemma}

\begin{proof}
Suppose that $f: \mathbb{P} \rightarrow \mathbb{P}$ is an automorphism.  Pullback by $f$ yields an isomorphism of class groups $f^*: \mathrm{Cl}(\mathbb{P}) \rightarrow \mathrm{Cl}(\mathbb{P})$.  Since $\mathrm{Cl}(\mathbb{P}) \cong \Z$ and ampleness of classes must be preserved, $f^* \mathcal{O}(1) = \mathcal{O}(1)$.  It follows that for every $m$ we have an isomorphism $f^*: H^0(\mathbb{P},\mathcal{O}(m)) \rightarrow H^0(\mathbb{P},\mathcal{O}(m))$ and that these isomorphisms are compatible with tensor product.  Since $H^0(\mathbb{P},\mathcal{O}(m))$ is the $m$th graded piece of $S$, these assemble to a graded isomorphism $f^*: S \rightarrow S$, which induces the original $f$.  It's straightforward to check that only maps of the form $x_i \mapsto t^{a_i}x_i$ for all $i$ are in the kernel, and that every $t \neq 1$ gives a non-identity map $S \rightarrow S$ by well-formedness.
\end{proof}

Given a fixed embedding of a weighted projective hypersurface $X \subset \mathbb{P}$, we'll define the subgroup $\mathrm{Lin}(X) \subset \mathrm{Aut}(X)$ of {\it linear automorphisms} to consist of those automorphisms of $X$ which extend to $\mathbb{P}$ (or, equivalently, to some graded automorphism of $S$).  We retain the terminology ``linear" by analogy with ordinary projective space, but it's important to note that the images of the variables $x_i$ under a graded automorphism of $S$ need not be linear as polynomials. For example, if $S = \C[x_0,x_1,x_2]$ and $x_0,x_1,x_2$ have weights $4,3,$ and $1$, respectively, we could define an element of $\mathrm{Aut}(S)$ by $x_0 \mapsto x_0 + x_1 x_2 - x_2^4$, $x_1 \mapsto 2x_1 + x_2^3$, $x_2 \mapsto - x_2$.

We'll need the following fact about automorphisms of graded polynomial rings in Sections \ref{bounds} and \ref{verygeneral}.  We say that a group $G \subset \mathrm{Aut}(S)$ is {\it diagonalizable} if, after conjugating $G$ by an automorphism of $S$, each element $g$ of $G$ sends each variable $x_i$ to some scalar multiple of itself.  Any such automorphism also sends an arbitrary monomial to a scalar multiple of itself.

\begin{lemma} \label{diagonalizable}
Let $S = \C[x_0,\ldots,x_{n+1}]$ be a weighted polynomial ring with weights $a_0, \ldots, a_{n+1}$ and $\mathrm{Aut}(S)$ the group of graded automorphisms.  If $G \subset \mathrm{Aut}(S)$ is a finite abelian group, then $G$ is diagonalizable.
\end{lemma}

This is, of course, a generalization of the well-known fact that any finite abelian group in $\mathrm{GL}_{n+2}(\C)$ is diagonalizable.

\begin{proof}
The group $\mathrm{Aut}(S)$ embeds naturally in $\bigoplus_{b \in \mathcal{B}} \mathrm{GL}(S_b)$, where $\mathcal{B} = \{b : b = a_i \text{ for some } i\}$ is the set of integers that appear as a weight of $S$.  (In particular, $\mathrm{Aut}(S)$ is a linear algebraic group.) To diagonalize a finite abelian group $G \subset \mathrm{Aut}(S)$, we'll focus on the representation on each of these pieces $S_b$ in turn.  Of course, within each $S_b$, we can diagonalize $G$ using some change of coordinates in $S_b$, but we must prove that we can do this simultaneously for all $b$ using conjugation by an element $\gamma$ in $\mathrm{Aut}(S)$.

To do this, we put the integers in $B$ in increasing order and construct $\gamma$ inductively.  For the smallest integer $b_0 \in \mathcal{B}$, there is a basis of $S_{b_0}$ given by $\{x_i: a_i = b_0 \}$.  Define $\gamma$ on the smallest weight variables in such a way that the representation of $G$ on $S_{b_0}$ becomes diagonal in the given basis (this is no problem because all monomials of degree $b_0$ are generators of $S$).  By the inductive hypothesis, we now assume $G$ acts by multiplication by a scalar on all variables $x_i$ of weight $a_i < b$ and consider the action on $S_b$.  Write a basis for $S_b$ beginning with the variables of weight $b$, followed by the other monomials of degree $b$.  By the inductive hypothesis, we've already constructed a $\gamma$ so that the representation of $G$ after changing coordinates in only the smaller variables is of the form
$$
\begin{pmatrix}
    A(g) & 0 \\
    B(g) & D(g)
\end{pmatrix}.$$
Here $D(g)$ is a diagonal sum of characters of $G$, because each $g \in G$ acts on monomials in the smaller weight variables by scalar multiplication.  Since the entire space $S_b$ must also be a direct sum of one-dimensional characters of $G$, we can find a change of coordinates affecting only coordinates in the first part of the basis (variables of weight $b$) so that the representation becomes diagonal.  This finishes the definition of the automorphism $\gamma$ on variables of weight up to $b$.  By induction, the proof is complete.
\end{proof}

\section{Linearity of Hypersurface Automorphisms} \label{linearity}

Let $X$ be a smooth hypersurface in $\mathbb{P}^{n+1}$.  Then, automorphisms of $X$ extend to $\mathbb{P}^{n+1}$ whenever $n \geq 1$ and $d \geq 3$ unless $(n,d) = (1,3)$ or $(2,4)$.  When $n \geq 3$, this is a consequence of the Grothendieck-Lefschetz theorem \cite[Expos\'{e} XII, Corollaire 3.6]{SGA2}, which holds even for arbitrarily singular hypersurfaces in projective space. For $n = 2$, $d \neq 4$, this is a theorem due to Matsumura and Monsky \cite{MM}; when $n = 1, d \geq 4$, it is due to Chang \cite{Chang}.  In this section, we prove a generalization to hypersurfaces in weighted projective space. We will deduce this statement from the local version of Grothendieck's Lefschetz theorems, but some care is needed. One complication is that a hypersurface in weighted projective space need not be a Cartier divisor (this is usually assumed even for variants of the global Lefschetz theorems that deal with singular varieties, e.g. \cite{RS}).  

\begin{theorem}
\label{main-theorem}
Let $X \subset \mathbb{P}(a_0,\ldots,a_{n+1})$ and $X' \subset \mathbb{P}(a_0',\ldots,a_{n+1}')$ be two complex weighted projective hypersurfaces of weighted degrees $d$ and $d'$, respectively.  Suppose further that $X$ and $X'$ are well-formed and quasismooth, neither is a linear cone, and one of the following holds:
\begin{enumerate}
    \item $n \geq 3$; or
    \item $n = 2$ and $a_0 + a_1 + a_2 + a_3 \neq d$.
\end{enumerate}
Then, if $g: X' \rightarrow X$ is an isomorphism, we have $d = d'$, the $a_i$ coincide with the $a_i'$ up to rearrangement, and $g$ is induced by an automorphism of $\mathbb{P}(a_0,\ldots,a_{n+1})$.
\end{theorem}

\begin{remark}
    \begin{enumerate}
        \item The assumption of well-formedness is necessary (note that a quasismooth hypersurface in a well-formed weighted projective space of dimension at least $3$ is automatically itself well-formed \cite[Theorem 6.17]{Iano-Fletcher}).  For example, any hypersurface $X_4 \subset \mathbb{P}^4(2,2,2,2,2)$ is isomorphic as a variety to a hypersurface $X_2 \subset \mathbb{P}^4(1,1,1,1,1)$ with the same equation.
        \item A hypersurface is a {\it linear cone} if the equation $f$ contains a linear term $x_i$ for some $i$. In this case, the hypersurface $f = 0$ is isomorphic to a weighted projective space of smaller dimension, and the theorem above fails rather trivially: for example, $\{x_0 = 0 \} \subset \mathbb{P}^4(1,1,1,1,1)$ and $\{x_0 = 0\} \subset \mathbb{P}^4(2,1,1,1,1)$ are both isomorphic to $\mathbb{P}^3$.  We avoid the linear cone case to exclude this pathology.  
       \item The part of the theorem which states that the degree and weights coincide for any two embeddings of the same weighted projective hypersurface was shown in the greater generality of quasismooth weighted complete intersections of dimension at least $3$ by Przyjalkowski and Shramov \cite[Proposition 1.5]{PrS}.  They also proved in another paper that the action of any \textit{linearly reductive} algebraic group $\Gamma$ on a quasismooth weighted complete intersection of dimension at least $3$ or of dimension $2$ with nontrivial canonical class comes from an action of $\Gamma$ on the ambient weighted projective space \cite[Theorem 1.3]{PrS2}.  We remove the reductivity assumption, answering \cite[Question 1.5]{PrS2} for hypersurfaces.  Some results on the class groups of weighted complete intersections of dimension at least $3$ that we'll use below date back further to \cite{Mori, BR}.
    \end{enumerate}
\end{remark}

\begin{proof}
Let $S = \C[x_0,\ldots,x_{n+1}]$ be the graded polynomial ring in variables $x_i$ of degrees $a_i$, and define $S'$ in the same way, so that $\mathbb{P} \coloneqq \mathbb{P}(a_0,\ldots,a_{n+1}) = \mathrm{Proj}(S)$ and $\mathbb{P}' \coloneqq \mathbb{P}(a_0',\ldots,a_{n+1}') = \mathrm{Proj}(S')$.  If $f$ and $f'$ are the homogeneous polynomials of weighted degrees $d$, $d'$ defining $X$ and $X'$, respectively, then $X = \on{Proj}(S/(f))$ and $X' = \on{Proj}(S'/(f'))$.  We'll first show that the isomorphism $g: X' \rightarrow X$ comes from an isomorphism of graded rings $g^*: S/(f) \rightarrow S'/(f')$. This is non-trivial since not every isomorphism of two Proj's comes from an underlying isomorphism of the graded rings.  We use the notation $\mathcal{O}_{\mathbb{P}}(m)$ and $\mathcal{O}_X(m)$ for the sheaves coming from each respective Proj construction and $i: X \rightarrow \mathbb{P}$ for the closed immersion coming from the quotient of graded rings $S \rightarrow S/(f)$.

\begin{lemma}
\label{proj_sheaves}
In the setting of \Cref{main-theorem}, for each $m \geq 0$, the natural map $i^* \mathcal{O}_{\mathbb{P}}(m) \rightarrow \mathcal{O}_X(m)$ of sheaves is an isomorphism.  Furthermore, the sheaf $\mathcal{O}_X(m)$ is reflexive on $X$.
\end{lemma}

\begin{proof}
We may check these conditions on each affine open chart $D^+_X(x_j) \subset X, j = 0, \ldots, n+1$ where $x_j$ doesn't vanish. (Here we note that each $x_j$ is still nonzero in the quotient $S/(f)$ because $X$ is not a linear cone.) These cover $X = \mathrm{Proj}(S/(f))$ and the corresponding collection of opens $D^+_{\mathbb{P}}(x_j)$ also covers $\mathbb{P}$.  Use the notation $S_{(\alpha)}$ for $\alpha \in S$ a homogeneous element of positive degree to mean the degree zero part of the localization $S_{\alpha}$; similarly, when $M$ is a graded $S$-module, $M_{(\alpha)}$ means the degree zero part of the localization.  We then have $D^+_{\mathbb{P}}(x_j) = \Spec(S_{(x_j)})$, $D^+_X(x_j) = \Spec((S/(f))_{(x_j)})$, and that $\mathcal{O}_{\mathbb{P}}(m)|_{D^+_{\mathbb{P}}(x_j)}$ is the quasicoherent sheaf associated to the $S_{(x_j)}$-module $(S(m))_{(x_j)}$.  The claimed isomorphism amounts to checking
$$(S(m))_{(x_j)} \otimes_{S_{(x_j)}} (S/f)_{(x_j)} \cong ((S/f)(m))_{(x_j)}.$$
This follows from the commutativity of localizations and tensor product, as well as keeping track of the degrees on each side.  To show that $\mathcal{O}_X(m)$ is reflexive, we use the following criterion \cite[Proposition 1.6]{Hartshorne}: a coherent sheaf $\mathcal{F}$ on a normal, integral scheme $X$ is reflexive iff 1) it is torsion-free, and 2) for every open set $U \subset X$ and closed $Y \subset U$ of codimension at least $2$, the restriction $\mathcal{F}(U) \rightarrow \mathcal{F}(U \setminus Y)$ is a bijection.

 It's clear that $\mathcal{O}_X(m)$ is torsion-free because $S/(f)$ is (otherwise, the affine cone would be non-reduced, contradicting quasismoothness). We can check the second property on the affine charts $D^+_X(x_j)$.   On this chart, $X$ is a quotient of the variety $X' \subset \mathbb{A}_{x_0,\ldots,\hat{x}_j,\ldots, x_{n+1}}^{n+1}$ by $\mu_{a_j}$, where $X'$ is cut out by the equation $f(x_0,\ldots,x_{j-1},1,x_{j+1},\ldots,x_{n+1}) = 0$.  Denote by $q:  X' \rightarrow X$ the quotient morphism.  By well-formedness of $X$, this action is free away from a codimension $2$ subset $Z \subset X'$.  Therefore, given $s \in \mathcal{O}_{D^+_X(x_j)}(U \setminus Y)$, the restriction of $s$ to $U \setminus (Y \cup q(Z))$ lifts to a regular function on the preimage of this set in $X'$.  By quasismoothness, $X'$ is smooth, and in particular normal, so this function extends to all of the preimage of $U$ and remains homogeneous of the same weight for the $\mu_{a_j}$-action.  Hence $s$ extends to $U$, completing the proof.
\end{proof}

We require the reflexivity of $\mathcal{O}_X(m)$ in order to view it as a member of the class group $\mathrm{Cl}(X)$ below.

For each integer $m$, we claim that the following sequence of sheaves on $\mathbb{P}$ is exact:
\begin{equation}
\label{tricky_exact_sequence}
    0 \rightarrow \mathcal{O}_{\mathbb{P}}(m-d) \xrightarrow{\cdot f} \mathcal{O}_{\mathbb{P}}(m) \rightarrow \mathcal{O}_X(m) \rightarrow 0.
\end{equation}
Here the second map is multiplication by $f$.  Indeed, the same sequence with $\mathcal{O}_X(m)$ replaced by $\mathcal{O}_\mathbb{P}(m)|_X$ is manifestly exact because it corresponds to the exact sequence of modules $0 \rightarrow S(m-d) \xrightarrow{\cdot f} S(m) \rightarrow (S/(f))(m) \rightarrow 0$.  By \Cref{proj_sheaves}, $\mathcal{O}_\mathbb{P}(m)|_X = i^*\mathcal{O}_{\mathbb{P}}(X) \cong \mathcal{O}_X(m)$, so \eqref{tricky_exact_sequence} is also exact.

By \cite[Section 1.4]{Dolgachev}, $H^1(\mathbb{P},\mathcal{O}_{\mathbb{P}}(m-d)) = 0$ for all $m$ because the dimension of the weighted projective space $\mathbb{P}$ is at least $2$. This is the analog of projective normality for weighted projective hypersurfaces. Therefore, we arrive at the exact sequence of global sections
$$0 \rightarrow H^0(\mathbb{P},\mathcal{O}_{\mathbb{P}}(m-d)) \xrightarrow{\cdot f} H^0(\mathbb{P},\mathcal{O}_{\mathbb{P}}(m)) \rightarrow H^0(X,\mathcal{O}_X(m)) \rightarrow 0.$$
Since $S_m \cong \mathcal{O}_{\mathbb{P}}(m)$ and the image of the first map is $f S_{m-d}$, we may identify $H^0(X,\mathcal{O}_X(m))$ with the $m$th graded piece of $S/(f)$.  Therefore,
$$S/(f) = \bigoplus_{m=0}^{\infty} H^0(X,\mathcal{O}_X(m)),$$
and similarly with $S'/(f')$. Thus, it will be enough to find maps $H^0(X,\mathcal{O}_X(m)) \rightarrow H^0(X',\mathcal{O}_{X'}(m))$ for each $m$ to build our desired map $S/(f) \rightarrow S'/(f')$.

Since $X$ is quasismooth, the punctured affine cone $C_X^* \subset \mathbb{A}^{n+2} \setminus \{0\}$ is smooth.  The hypersurface $X$ is a quotient of $C_X^*$ by the action of $G \coloneqq \mathbb{G}_{\on{m}}$ on $\mathbb{A}^{n+2} \setminus \{0\}$; we'll denote the quotient morphism by $q: C_X^* \rightarrow X$.  Using the assumption of well-formedness of $X$, $\on{Sing}(X)$ is codimension at least $2$ in $X$.  $\on{Sing}(X)$ is also the image of the locus in $C_X^*$ where the $G$-action has nontrivial stabilizers.

Thus, we have $\on{Cl}(X) \cong \on{Cl}(X_{\on{sm}})$, where $X_{\on{sm}} = X \setminus \on{Sing}(X)$ is the smooth locus.  Because $X_{\on{sm}}$ is smooth, $\on{Cl}(X_{\on{sm}}) = \on{Pic}(X_{\on{sm}})$.  Further, $\on{Pic}(X_{\on{sm}})$ is isomorphic to the group $\on{Pic}^G(q^{-1}(X_{\on{sm}}))$ of $G$-equivariant vector bundles on $q^{-1}(X_{\on{sm}})$, since $q^{-1}(X_{\on{sm}}) \rightarrow X_{\on{sm}}$ is the quotient by a free group action.  Finally, $q^{-1}(X_{\on{sm}})$ has complement of codimension at least $2$ in $C_X^*$, so $G$-equivariant line bundles on $X_{\on{sm}}$ extend to all of $C_X^*$.  In summary,
$$\on{Cl}(X) \cong \on{Pic}(X_{\on{sm}}) \cong \on{Pic}^G(q^{-1}(X_{\on{sm}})) \cong \on{Pic}^G(C_X^*).$$
The isomorphism $g : X' \rightarrow X$ induces an isomorphism $X'_{\on{sm}} \rightarrow X_{\on{sm}}$, which also gives a pullback map $\on{Pic}(X_{\on{sm}}) \xrightarrow{\cong} \on{Pic}(X'_{\on{sm}})$.  We'll identify $\mathrm{Cl}(X)$ and $\mathrm{Pic}^G(C_X^*)$ below without further comment and use the same notation $\mathcal{O}_X(m)$ for the corresponding elements in either of these groups.

\begin{proposition}
\label{preservegenerator}
Suppose that $X$ and $X'$ satisfy the conditions of \Cref{main-theorem}.  Then an isomorphism $g: X' \rightarrow X$ induces an isomorphism $\on{Cl}(X) \xrightarrow{\cong} \on{Cl}(X')$ which maps $\mathcal{O}_X(1)$ to $\mathcal{O}_{X'}(1)$.
\end{proposition}

\begin{proof}
We argue differently depending on the dimension of $X$. When $\dim(X) \geq 3$, we claim that $\on{Pic}^G(C_X^*) \cong \Z \cdot \mathcal{O}_X(1)$ (and analogously for $X'$).

Indeed, if we forget the $G$-equivariant structure, all line bundles on the smooth variety $C_X^*$ are trivial when $\dim(X) \geq 3$.  This is because the local ring $\mathcal{O}_{C_X,0}$ is a complete intersection ring of dimension at least $4$, regular outside its maximal ideal, so it is a parafactorial ring \cite[Expos\'{e} XI, Th\'{e}or\`{e}me 3.13]{SGA2}.  It follows that $\on{Pic}(C_X) \rightarrow \on{Pic}(C_X^*)$ is an isomorphism; furthermore, $\on{Pic}(C_X)$ is trivial \cite[Section 3.2.2]{Dolgachev}. From $\mathrm{Pic}(C_X^*) = 0$, we deduce that the group of $G$-equivariant line bundles on $C_X^*$ is naturally isomorphic to the character group of $G$, namely $\Z$ \cite[Lemma 4.1.7]{Brion}.  It's straightforward to check that the linearization of the trivial bundle by the character $t \mapsto t^m$ coincides with the $G$-equivariant bundle $\mathcal{O}_X(m)$.  Therefore, $\on{Cl}(X) = \Z \cdot \mathcal{O}_X(1)$.  Since $X$ and $X'$ are isomorphic, $\dim(X') \geq 3$ also, so identical reasoning shows $\on{Cl}(X') = \Z \cdot \mathcal{O}_{X'}(1)$.

The pullback of an ample divisor by an isomorphism is still ample, so the isomorphism $\on{Cl}(X) \xrightarrow{\cong} \on{Cl}(X')$ must send $\mathcal{O}_X(1)$ to $\mathcal{O}_{X'}(1)$.  This proves the proposition when $n \geq 3$.

If $\dim(X) = 2$, we use a different argument, since we may not have $\on{Cl}(X) \cong \Z$ in this case.  The canonical class of $X$ is $K_X = \mathcal{O}_X(r)$, where $r = d - \sum_i a_i \neq 0$ by assumption.  We also have $K_{X'} = \mathcal{O}_X(r')$ with $r' = d' - \sum_i a_i'$ having the same sign as $r$, depending on whether $X$ and $X'$ have ample or anti-ample canonical class.  The canonical class is preserved by isomorphism, so pullback sends $\mathcal{O}_X(r)$ to $\mathcal{O}_{X'}(r')$.  We'll prove that $r = r'$ and moreover that $\mathcal{O}_X(1)$ maps to $\mathcal{O}_{X'}(1)$.

\begin{lemma}
\label{fundamentalgroup}
Let $V$ be a connected scheme of finite type over $\C$.  If $\pi_1(V) = 1$, then $\Pic(V)$ is torsion-free.
\end{lemma}

\begin{proof}
For any positive integer $\ell$, we have the following Kummer exact sequence of sheaves of abelian groups on $V$ in the \'{e}tale topology:
$$1 \rightarrow \mu_{\ell} \rightarrow \mathbb{G}_{\on{m},V} \xrightarrow{x \mapsto x^{\ell}} \mathbb{G}_{\on{m},V} \rightarrow 1.$$
The associated long exact sequence in cohomology gives
$$\cdots \rightarrow H^1(V,\mu_{\ell}) \rightarrow \on{Pic}(V) \xrightarrow{L \mapsto L^{\ell}} \on{Pic}(V) \rightarrow \cdots.$$
Since $V$ is connected, $H_1(V,\Z)$ is the abelianization of $\pi_1(V)$; hence $H_1(V,\Z) = 0$. The universal coefficient theorem then gives that $H^1(V,\mu_{\ell}) \cong H^1(V,\Z/{\ell}) = 0$, so that the $\ell$th tensor power map on $\on{Pic}(V)$ is injective.  Since this holds for any positive integer $\ell$, $\on{Pic}(V)$ is torsion-free.
\end{proof}

\begin{commentA}
Here's another, more geometric proof of \Cref{fundamentalgroup}.  We use the cyclic covering trick.  Let $L$ be a torsion element of $\Pic(V)$ of order $\ell$, so that there is an isomorphism $s: L^{\ell} \cong \mathcal{O}_V$ given by some nowhere vanishing global section of $L^{\ell}$.  We define a finite covering
$$\underline{\mathrm{Spec}}_V \bigoplus_{i = 0}^{\ell-1} L^{i} \rightarrow V.$$
Multiplication in the sheaf of $\mathcal{O}_V$-algebras $\mathcal{A} \coloneqq \oplus_{i = 0}^{\ell-1} L^{i}$ is defined by the usual multiplication $L^{a} \otimes L^{b} \rightarrow L^{\otimes (a+b)}$ when $a + b < \ell$ and by the composition $L^a \otimes L^b \rightarrow L^{a+b} \xrightarrow{s} L^{a+b-\ell}$ if $a + b \geq \ell$.  Then this is an \'etale cover of $V$.  However, since the fundamental group of $V$ is trivial, this covering must simply be a disjoint union of copies of $V$.  This means that $L \cong \mathcal{O}_V$, as required.
\end{commentA}

Now we'll conclude the proof of \Cref{preservegenerator}. By \cite[Section 3.2.2]{Dolgachev}, the fundamental group of the punctured affine cone $C_X^*$ vanishes whenever $X$ is a quasismooth hypersurface of dimension at least $2$.  Therefore, the (non-equivariant) Picard group $\on{Pic}(C_X^*)$ is torsion-free.  Next, we use the exact sequence
$$0 \rightarrow \Z \rightarrow \on{Pic}^G(C_X^*) \rightarrow \on{Pic}(C_X^*),$$
where the first map sends $1 \mapsto \mathcal{O}_X(1)$ and the second map forgets the linearization.  Exactness of this sequence follows from \cite[Theorem 4.2.2]{Brion} plus the observation that $\Z \rightarrow \on{Pic}^G(C_X^*)$ is injective by ampleness of $\mathcal{O}_X(1)$ on $X$.  It follows that $\on{Pic}^G(C_X^*)$ is torsion-free, since any torsion element must map to zero in $\on{Pic}(C_X^*)$, and hence be in the image of $\Z \rightarrow \on{Pic}^G(C_X^*)$.  In addition, we see that $\mathcal{O}_X(k) \in \on{Pic}^G(C_X^*)$ is not divisible by any integer other than the factors of $k$, because otherwise the cokernel of the first map would contain torsion elements.

This implies $r = r'$ above.  Indeed, if $g^*(\mathcal{O}_X(r)) \cong \mathcal{O}_{X'}(r')$, then the image $L \in \on{Pic}^G(C_{X'}^*)$ of $\mathcal{O}_X(1)$ satisfies $L^r = \mathcal{O}_{X'}(r')$.  Since $\mathcal{O}_{X'}(r')$ is divisible by $r$, $r$ is a factor of $r'$ by the above.  Arguing symmetrically with the inverse isomorphism gives $r = r'$.  In particular, we then have that the image of $\mathcal{O}_X(1)$ differs from $\mathcal{O}_{X'}(1)$ by a torsion element of order $r$.  Because we also saw that $\on{Pic}^G(C_X^*)$ is torsion free, this proves that $g^*(\mathcal{O}_X(1)) \cong \mathcal{O}_{X'}(1)$. \end{proof}

Finally, we're ready to complete the proof of \Cref{main-theorem}.  We've now seen that in the assumptions of the theorem, pullback by the isomorphism $g$ sends $\mathcal{O}_X(m)$ to $\mathcal{O}_{X'}(m)$ for any $m \in \Z$.  Therefore, $g$ induces isomorphisms
$$g^*: H^0(X,\mathcal{O}_{X}(m)) \rightarrow H^0(X',\mathcal{O}_{X'}(m)).$$
These maps respect the tensor product of sections, so we may assemble them into an isomorphism of graded rings $g^*: S/(f) \rightarrow S'/(f')$.  In particular, $g^*$ gives an isomorphism of the maximal irrelevant ideal $\mathfrak{m}$ of elements of positive degree in $S/(f)$ with $\mathfrak{m}' \subset S'/(f')$.  The number of generators of $\mathfrak{m}$ and their degrees (up to reordering) coincide with those of $\mathfrak{m'}$.  This is because if $m$ is the smallest positive index with the graded piece $\mathfrak{m}_m$ nonempty, the number of generators of degree $m$ equals $\dim(\mathfrak{m}_m)$, and then we can factor out by these generators and repeat inductively.

Since $X$ is quasismooth and not a linear cone, we have that $d$, the degree of the smallest relation among the $x_i$, is strictly greater than all weights $a_i$.  It follows that $x_0,\ldots,x_{n+1}$ are a minimal system of $n+1$ generators for the homogeneous ideal $\mathfrak{m}$.  Similarly, $x_0',\ldots,x_{n+1}'$ generate $\mathfrak{m}'$, which is an isomorphic to $\mathfrak{m}$, so this set of generators must also be minimal and the collection $a_0,\ldots,a_{n+1}$ must be the same as $a_0',\ldots,a_{n+1}'$, up to reordering.  Dimension counting yields that the relations are then also in the same degree, so we have $d = d'$.

Lastly, we observe that the isomorphism $g^* \colon S/(f) \rightarrow S'/(f')$ is induced by an isomorphism $S \rightarrow S' $.  Indeed, for $m < d$, the $m$th graded piece of $S$ is isomorphic to that of $S/(f)$, so $g^*$ gives isomorphisms $S_m \rightarrow S'_m$.  All the generators of $S$ are contained in $S_m$ for $m < d$ so this gives rise to a homomorphism $S \rightarrow S'$.  The inverse of $g$ similarly gives a morphism $S' \rightarrow S$ whose composition with the previous map is the identity on generators, and hence on all of $S$.  Therefore, the homomorphism $S \rightarrow S'$ is an isomorphism.  By our work above, both $\on{Proj}(S)$ and $\on{Proj}(S')$ are isomorphic to $\mathbb{P}(a_0,\ldots,a_{n+1})$, so the isomorphism $S \rightarrow S'$ gives an automorphism of this weighted projective space inducing the original isomorphism $g \colon X' \rightarrow X$ of hypersurfaces, as claimed.
\end{proof}

This theorem fails if we weaken the hypotheses on dimension.  As mentioned above, two smooth plane curves $C$ and $C'$ in $\mathbb{P}^2$ of degree at least $4$ are isomorphic if and only if they differ by an isomorphism of the projective plane \cite{Chang}.  However, the situation for weighted projective curves is considerably more complicated: there exist many curves of genus at least $2$ which can be embedded as well-formed quasismooth hypersurfaces in different weighted projective spaces.

\begin{example}[Hyperelliptic Curves]
Let $C$ be a hyperelliptic curve of genus $g$, $p: C \rightarrow \mathbb{P}^1$ a $2:1$ cover of $\mathbb{P}^1$, and $P \in C$ a ramification point of the cover.  Then $C$ is isomorphic to $\mathrm{Proj} (R(C,P))$, where
$$R(C,P) \coloneqq \bigoplus_{k=0}^{\infty} H^0(C,kP).$$

This ring has generators $x,y$, and $z$ in degrees $1$, $2$, and $2g+1$, respectively, and a single relation in degree $4g+2$.  It is possible to choose the generator $z$ so that the relation has the form $f(x,y,z) \coloneqq z^2 - h(x^2,y) = 0$, where $h$ is a polynomial of degree $2g+1$.  This embeds $C$ as a quasismooth hypersurface of degree $4g+2$ in $\mathbb{P}^2(2g+1,2,1)$.  However, if we use $R(2C,P) = R(C,K_C)$ instead, we have another embedding of the same curve as a degree $2g+2$ hypersurface in $\mathbb{P}^2(g+1,1,1)$.
\end{example}

\begin{example}
There are also non-hyperelliptic curves exhibiting similar behavior.  For $C$ a smooth non-hyperelliptic curve of genus $3$, we have the canonical embedding of $C$ in $\mathbb{P}^2 = \mathbb{P}^2(1,1,1)$ as a degree $4$ plane curve.  Suppose further that $C$ has the property that there is a line $\ell \subset \mathbb{P}^2$ tangent to $C$ at a point $P$ with multiplicity $4$.  Then we have $4P \sim K_C$.  One can show that ring $R(C,P)$ has generators in degrees $1$, $3$, and $4$, giving an embedding $C \subset \mathbb{P}^2(4,3,1)$ as a hypersurface of degree $12$.  There are many similar examples for curves of higher genus.
\end{example}

When $\dim X = 2$ but $\sum_i a_i = d$, \Cref{main-theorem} also fails.  This is because the resulting hypersurfaces are K3 surfaces in this case, which frequently have infinite automorphism group. An example of Fano and Severi of a quartic surface in $\mathbb{P}^3$ with infinite automorphism group is described in the proof of Theorem 4 in \cite{MM}, for instance.  However, we'll see below in \Cref{finite} that the linear automorphism groups of weighted projective surfaces with $\sum_i a_i = d$ are always finite.  Hence, some automorphisms of these surfaces aren't linear.

\section{Bounds on Linear Automorphism Groups} \label{bounds}

Recall that the linear automorphism group $\on{Lin}(X) \subset \on{Aut}(X)$ of a hypersurface $X \subset \mathbb{P}$ is the subgroup of automorphisms that extend to $\mathbb{P}$. As long as the degree of $X$ exceeds all weights of $\mathbb{P}$, the only automorphism of $\mathbb{P}$ that fixes $X$ pointwise is the identity (use the same argument on extending morphisms of graded rings as above), so in this case we may consider $\on{Lin}(X)$ as a subgroup of both $\mathrm{Aut}(X)$ and $\mathrm{Aut}(\mathbb{P})$. \Cref{main-theorem} shows that $\on{Lin}(X) = \mathrm{Aut}(X)$ whenever $\dim(X) \geq 3$ or $\dim(X) = 2$ and $K_X \ncong \mathcal{O}_X$.

In this section, we prove several results on the size of $\on{Lin}(X)$. These will imply the same results for $\mathrm{Aut}(X)$ in most dimensions in light of the last section.  First, we give a criterion in terms of the degree $d$ and the weights $a_0,\ldots,a_{n+1}$ which determines whether or not $\on{Lin}(X)$ is finite.

\begin{theorem}
\label{finite}
Let $X \subset \mathbb{P}(a_0,\ldots,a_{n+1})$ be a well-formed, quasismooth weighted projective hypersurface of degree $d$, where $n \geq 1$.  Then $\on{Lin}(X)$ is finite if and only if one of the following two conditions holds:
\begin{enumerate}
    \item $d > 2 \max \{a_0,\ldots,a_{n+1}\}$; or
    \item $d = 2 \max \{a_0,\ldots,a_{n+1}\}$ and only a single weight achieves the maximum.
\end{enumerate}
Further, if neither (1) nor (2) holds (so that $\mathrm{Lin}(X)$ is infinite), then $X$ is rational.
\end{theorem}

\begin{remark}
\begin{enumerate}
    \item This generalizes a theorem of Matsumura and Monsky \cite[Theorem 1]{MM}, which shows that the linear automorphism group of a smooth hypersurface of degree $d$ in $\mathbb{P}^{n+1}$ over an algebraically closed field $k$ is finite if $n \geq 2$ and $d \geq 3$.  For $k = \C$, this result was known in some form at least as far back as 1880, when it appeared in a work of Jordan \cite{Jordan2} (see also \cite[Section 6]{OS} for further historical remarks on this theorem). Note also that if $d < 3$, $X$ is a hyperplane or a quadric.  These always have infinite linear automorphism groups and are rational for $n \geq 1$.
    \item  Some partial results along these lines were known previously.  In the special case when $X$ is a {\it smooth} (rather than just quasismooth) well-formed weighted projective hypersurface, $\dim(X) \geq 3$ and $K_X \ncong \mathcal{O}_X$, it was proven in \cite[Corollary 1.4]{PrS} that $\mathrm{Aut}(X)$ is finite unless $X$ is isomorphic to either $\mathbb{P}^n$ or a quadric hypersurface of $\mathbb{P}^{n+1}$.  Note that $X$ is never smooth unless the weights are pairwise relatively prime and all divide the degree.  However, the result of \cite{PrS} also applies to smooth weighted complete intersections with these properties; their methods are rather different than ours.  A statement similar to the first part of \Cref{finite} also appears in work of Bunnett \cite[Theorem 3.13]{Bunnett}, but gives an incorrect degree bound of $d \geq \max \{a_0,\ldots,a_{n+1}\} + 2$ for finiteness.
\end{enumerate}
\end{remark}

\begin{proof}[Proof of \Cref{finite}]
First, assume that one of the two conditions on $d$ in \Cref{finite} holds.  We'll show that $\on{Lin}(X)$ is finite, using generally the same approach as in \cite{MM}.  By \Cref{autP}, any automorphism of $\mathbb{P}$ comes from a graded automorphism of $S = \C[x_0,\ldots,x_{n+1}]$, so we perform most of our analysis on the level of graded ring automorphisms.  For any graded homomorphism $S \rightarrow S$, the image of each generator $x_i$ is contained in the finite-dimensional vector space $S_{a_i}$.  Thus, the endomorphism monoid of $S$ is isomorphic to $\mathbb{A}^{N}$ as a variety, where $N \coloneqq \sum_i \dim(S_{a_i})$.  The linear algebraic group $\mathrm{Aut}(S)$ is an open subvariety of $\mathbb{A}^N$. This is a generalization of the fact that $\mathrm{GL}_{n+2}(\C)$ is an open subvariety of $\mathbb{A}^{(n+2)^2}$.

We saw in \Cref{autP} that $\mathrm{Aut}(S)$ acts on $\mathbb{P}(a_0,\ldots,a_{n+1})$ with kernel isomorphic to $\C^*$, where $t \in \C^*$ acts as $t \cdot x_i = t^{a_i} x_i$ for each $i$.  Let $G \subset \mathrm{Aut}(S)$ be the subgroup of elements mapping the polynomial $f$ defining $X$ to a multiple of itself.  Then $\on{Lin}(X) = G/\C^*$. Since $\on{Lin}(X)$ is an algebraic group, if it has trivial Lie algebra, then it must be finite.  We'll show that the Lie algebra $\mathfrak{g}$ of $G$ equals that of the subgroup $\ker(\mathrm{Aut}(S) \rightarrow \mathrm{Aut}(\mathbb{P})) \cong \C^*$; this implies that the Lie algebra of the quotient $\mathrm{Lin}(X) = G/\C^*$ is trivial, as required.

The tangent space to $\mathrm{Aut}(S)$ at the identity is naturally isomorphic to $S_{a_0} \oplus \cdots \oplus S_{a_{n+1}}$, where an element $z \coloneqq (z_0,\ldots,z_{n+1})$ corresponds to the infinitesimal automorphism $x_i \mapsto x_i + \epsilon z_i$.  Our aim is to show that if $z \in \mathfrak{g}$, then in fact $z$ is a multiple of $(a_0 x_0, \ldots, a_{n+1}x_{n+1})$, which is the derivative of the function $\C^* \rightarrow G$ given by $t \mapsto (x_i \mapsto t^{a_i} x_i)$ at $t = 1$.

Every $z$ in the Lie algebra of $\mathrm{Aut}(S)$ defines a derivation $D_z: S \rightarrow S$ by the formula $h \mapsto \frac{d}{d \epsilon} h(x + \epsilon z)|_{\epsilon = 0}$.  The fact that $z \in \mathfrak{g}$ means that $D_z(f) = c f$ for some constant $c$.  But we may express $D_z(f)$ in terms of partial derivatives $f_i \coloneqq \frac{\partial f}{\partial x_i}$ as:
$$D_Z(f) = \sum_i f_i z_i.$$
Therefore, we have
$$\sum_i f_i z_i = cf = \sum_i \frac{ca_i x_i}{d}f_i.$$
In this equation, the last equality follows from the following weighted variant of Euler's formula for homogeneous polynomials. Namely, for $f$ homogeneous of weighted degree $d$ in variables $x_i$ with weights $a_i$:
$$\sum_i a_i x_i f_i = df.$$
Rearranging the equation above now gives
\begin{equation}
\label{lincombo}
\sum_i \left(z_i -  \frac{ca_i x_i}{d} \right)f_i = 0.
\end{equation}
Since $X$ is a quasismooth hypersurface, its punctured affine cone in $\mathbb{A}^{n+2} \setminus \{0\}$ is smooth, so that the only common zero of the partial derivatives $f_0,\ldots,f_{n+1}$ is at the origin $x = 0$ in $\mathbb{A}^{n+2}$.  The ring $S$ is Cohen-Macaulay, each $f_i$ is a homogeneous element of positive degree in a graded ring, and these $n+2$ polynomials cut out the subvariety $\{0\}$ of codimension $n+2$ in $\mathbb{A}^{n+2}$. Therefore, it follows that any permutation of the $f_i$ form a regular sequence in the ring $S$.

A particular consequence of that fact is that $f_i$ is not a zero divisor in the ring $S/I_i$, where $I_i \coloneqq (f_0,\ldots,\hat{f}_i,\ldots,f_{n+1})$.  The equation \eqref{lincombo} then implies that $z_i - ca_i x_i/d \in I_i$ for each $i$.  That element is homogeneous of weighted degree $a_i$, so we can guarantee that it is zero if every nonzero polynomial in $I_i$ has degree greater than $a_i$.  Each partial derivative $f_j$ has degree $d - a_j$, so we can conclude that $z_i - ca_i x_i/d = 0$ if $d-a_j$ is greater than $a_i$ for all $j \neq i$.  Either of the two conditions in the theorem guarantees that this criterion is met for all $i$.  Therefore, if we assume one of these conditions, then $z_i = ca_i x_i/d = (c/d) a_i x_i$ for all $i$, as required.

Next, we'll show the converse: if either $d < 2 \max \{a_0,\ldots,a_{n+1}\}$ or $d = 2 \max \{a_0,\ldots,a_{n+1}\}$ and there are multiple weights equal to the maximum, then $\on{Lin}(X)$ is infinite.  Furthermore, we'll prove that $X$ is rational.  It's helpful to consider a few distinct cases.  Throughout, we'll assume the weights are arranged in decreasing order, so in particular $a_0 = \max \{a_0,\ldots,a_{n+1}\}$.  We note that either assumption on degree above guarantees that $d < a_0 + \cdots + a_{n+1}$ so that the hypersurface is Fano.

\begin{commentA}
Here's a full proof of the fact that hypersurfaces failing (1) and (2) in \Cref{finite} are Fano.

\begin{proposition}
    Let $X = \{f = 0\} \subset \mathbb{P}(a_0,\ldots,a_{n+1})$ be a quasismooth hypersurface of dimension at least $1$ which does not satisfy either of the conditions (1), (2) of \Cref{finite}.  Then $X$ is Fano.
\end{proposition}

\begin{proof}
Suppose first that $d < 2 \max \{a_0,\ldots,a_n\}$.  For quasismooth $X$, the minimum possible value of $d$ is $\max \{a_0,\ldots,a_n\}$, which corresponds to the case of a linear cone.  The hypersurface $X$ is then isomorphic to a weighted projective space, and is hence Fano.  If $d$ doesn't equal this minimum, suppose without loss of generality that $a_0$ is a maximum weight.  Then quasismoothness of $X$ means (by \Cref{monomialexistence}) that there must be a monomial $x_0 x_i$ of degree $d$ for some $i$.  Thus, $a_0 + a_i = d$.  Since $n \geq 1$, $d - \sum_i a_i < d- a_0 - a_i = 0$, so $X$ is Fano.  Similarly, if $d = 2 \max \{a_0,\ldots,a_n\}$, condition (2) of \Cref{finite} fails only if we have (at least) two weights $a_0$ and $a_1$ which equal $d/2$. In this case $d - \sum_i a_i < d - a_0 - a_1 < 0$, so $X$ is again Fano. 
\end{proof}
\end{commentA}

If $d = \max \{a_0,\ldots,a_{n+1}\} = a_0$, then $X$ is a linear cone, and hence isomorphic to $\mathbb{P}(a_1,\ldots,a_{n+1})$, which is rational.  We may also assume that the equation of $X$ is $x_0 + f(x_1,\ldots,x_n) = 0$, where $f$ is homogenenous of degree $a_0$.  Under the automorphism $x_0 \mapsto x_0 - f$ of $\mathbb{P}$, this becomes $x_0 = 0$.  Every automorphism of $\{x_0 = 0\} = \mathbb{P}(a_1,\ldots,a_{n+1})$ extends to $\mathbb{P}(a_0,\ldots,a_{n+1})$ and this group is infinite since $n \geq 1$.

Now suppose $d < \max \{a_0,\ldots,a_{n+1}\} < 2d$.  In order for $X$ to be quasismooth, its equation must contain a monomial $x_0 x_l$ (for some $l \neq 0$) with nonzero coefficient by \Cref{monomialexistence}.  By a transformation of $x_l$, we may assume that $x_0 x_l$ is the only term involving $x_0$. The equation then looks like 
$$x_0 x_l + x_l f(x_1,\ldots,x_{n+1}) + g(x_1,\ldots,\hat{x}_l,\ldots,x_n),$$
where $f$ is homogeneous of degree $d - a_l = a_0$ and may include $x_l$, while $g$ is homogeneous of degree $d$ and consists of terms not containing $x_0$ or $x_l$.  Composing with $x_0 \mapsto x_0 - f$ then eliminates the middle term.  After these transformations, it's clear that $X$ contains an infinite group of automorphisms
of the form $x_{0} \mapsto t x_{0}$,
$x_{l} \mapsto \frac{1}{t}x_{l}$, for $t \in {\mathbb C}^{*}$ (fixing all other variables).  To show that $X$ is rational, note that the $g$ term above is nonzero since $X$ is quasismooth.  On the open set $x_l \neq 0$, we may isolate $x_0 = -g/x_l$, so that the projection forgetting $x_0$ is a birational map to the rational toric variety $\mathbb{P}(a_1,\ldots,a_{n+1})$.

Finally, suppose $d = 2 \max \{a_0,\ldots,a_{n+1}\}$, but that both $a_0$ and $a_1$ equal $\frac{d}{2}$.  A similar series of reductions to the equation of $f$ works here: we can change variables so that the quadratic in $x_0$ and $x_1$ equals $x_0 x_1$ and eliminate any other terms involving $x_0$ and $x_1$.  The equation $x_0x_1 + f(x_2,\ldots,x_{n+1}) = 0$ has the same infinite family $x_0 \mapsto t x_0$, $x_1 \mapsto \frac{1}{t}x_1$ of automorphisms, and the projection forgetting the first coordinate is again a birational map.  This completes the proof.
\end{proof}

Next, we consider bounds on the size of the linear automorphism groups of hypersurfaces when they are finite.  Some results in this direction are known for hypersurfaces of degree $d$ in ordinary projective space $\mathbb{P}^{n+1}$, which we know have finite linear automorphism groups when $d \geq 3$.  An unpublished work of Bott and Tate from around 1961 showed that there is a bound on the size of the $\mathrm{Lin}(X)$ in terms of $d$ and $n$ (see \cite{OS} for an exposition of these ideas).  Around twenty years later, Howard and Sommese \cite{HS} showed that there is a constant $k_n$ for every dimension $n$ such that $|\mathrm{Lin}(X)| < k_n d^{n+1}$, for any $d \geq 3$.  We'll prove an even stronger theorem in the setting of weighted projective hypersurfaces.

\begin{theorem}
\label{bound}
For every positive integer $n$, there exists a constant $C_n$ depending only on $n$ with the following property: for any well-formed, quasismooth hypersurface $X \subset \mathbb{P}(a_0,\ldots,a_{n+1})$ of degree $d$ and dimension $n$, if $\on{Lin}(X)$ is finite, then
\begin{equation} \label{bound_eq}
    |\on{Lin}(X)| \leq C_n \frac{d^{n+1}}{a_0 \cdots a_{n+1}}.
\end{equation}
\end{theorem}

In particular, the same constant $C_n$ works for hypersurfaces in {\it any} weighted projective space of dimension $n$.  The comments following the proof of \Cref{Jordan_const} describe an explicit procedure for effectively producing a constant $C_n$ for which the theorem holds. We'll need the following definitions during the proof.

\begin{definition} \label{Jordan_def}
Let $\mathcal{G}$ be a group.  We say that $\mathcal{G}$ has the {\it Jordan property} if there exists a constant $J(\mathcal{G})$ such that for every finite subgroup $H \subset \mathcal{G}$, there exists a normal abelian subgroup $A \subset H$ with index $[H:A] \leq J(\mathcal{G})$.  The minimum $J(\mathcal{G})$ with this property is called the {\it Jordan constant} of $\mathcal{G}$.  The {\it weak Jordan constant} $\bar{J}(\mathcal{G})$ of $\mathcal{G}$ is the minimum constant such that every finite $H \subset \mathcal{G}$ has a (not necessarily normal) abelian subgroup $A \subset H$ with $[H:A] \leq \bar{J}(\mathcal{G})$.
\end{definition}

An 1878 result by Jordan \cite{Jordan} shows that $\mathrm{GL}_N(\C)$ has the Jordan property for all $N$ (for a modern exposition of his original proof and subsequent developments, see \cite{Breuillard}).  The explicit values of the Jordan constants $J_N \coloneqq J(\mathrm{GL}_N(\C))$ were not computed until much later: Collins \cite{Collins} calculated all the $J_N$ and in particular showed that when $N \geq 71$, $J_N = (N+1)!$; this index is achieved by the $N$-dimensional standard representation of the symmetric group $S_{N+1}$.  His proof relies on the classification of finite simple groups.

Weak Jordan constants have not been as well studied, but it follows from a theorem of A. Chermak and A. Delgado that for any group with the Jordan property, $\bar{J}(\mathcal{G}) \leq J(\mathcal{G}) \leq \bar{J}(\mathcal{G})^2$ (see \cite[Theorem 1.41]{Isaacs} and \cite[Remark 1.2.2]{PS}).  The precise values of $\bar{J}_N \coloneqq \bar{J}(\mathrm{GL}_N(\C))$ for small $N$ are computed in \cite{PS}, but to the author's knowledge there has been no complete calculation of all the $\bar{J}_N$.

\begin{proof}[Proof of \Cref{bound}]
We'll prove the theorem in the following three steps:
\begin{enumerate}
    \item[{\bf Step 1}:] Show that $\mathrm{Lin}(X)$ is the image of a finite group of graded ring automorphisms which fixes the function $f$ defining $X$ and has order $d|\mathrm{Lin}(X)|$.
    \item[{\bf Step 2}:] Find a uniform bound $C_n$ on the weak Jordan constants $\bar{J}(\mathrm{Aut}(S))$ of the graded automorphism groups of weighted polynomial rings $S$ in $n+2$ variables.
    \item[{\bf Step 3}:] Show that the order of an abelian group of graded ring automorphisms fixing $f$ is at most $d^{n+2}/(a_0 \cdots a_{n+1})$.
\end{enumerate}

\noindent {\bf Step 1}: Suppose that $G \coloneqq \on{Lin}(X)$ is a finite group, for a quasismooth hypersurface $X$ of degree $d$ in $\mathbb{P} = \mathbb{P}(a_0,\ldots,a_{n+1})$.  Let $S = \C[x_0,\ldots,x_{n+1}]$ be the weighted polynomial ring with weights $a_i$, $f$ be the polynomial defining the hypersurface $X$, and $\pi$ be the natural quotient homomorphism $\pi: \mathrm{Aut}(S) \rightarrow \mathrm{Aut}(\mathbb{P})$.  If $g \in \pi^{-1}(G) \subset \mathrm{Aut}(S)$, the induced automorphism of $\mathbb{P}$ preserves $f$, so $g \cdot f = cf$ for some constant $c \in \C$.  Define $H$ to be the subgroup of $\pi^{-1}(G)$ of elements $g$ which satisfy the stronger condition $g \cdot f = f$.  We claim that $\pi|_H: H \rightarrow G$ is a surjective homomorphism with kernel of order $d$, so that in particular $H$ is a finite group with $|H| = d |G|$.

It's clear that $\pi|_H: H \rightarrow G$ is surjective because we can compose any automorphism in $\pi^{-1}(G)$ with an element of $\ker(\pi) \cong \mathbb{C}^*$ (see \Cref{autP}) to scale the factor $c$ to $1$.  An element of the kernel of $\pi|_H$ is a $t \in \C^*$ with $f(t^{a_0}x_0,\ldots,t^{a_{n+1}}x_{n+1}) = t^d f = f$, so $t$ is a $d$th root of unity.  This proves $|H| = d |G|$. To bound the order of $G$, we can therefore analyze the group $H \subset \mathrm{Aut}(S)$ instead.

\noindent {\bf Step 2}:  
Next, we reduce to only considering abelian groups by computing the weak Jordan constant for the group $\mathrm{Aut}(S)$ of graded automorphisms. For any weighted polynomial ring $S$, $\mathrm{Aut}(S)$ is a linear algebraic group.  This implies that $\mathrm{Aut}(S)$ has the Jordan property.  However, even for a fixed number of variables, $\mathrm{Aut}(S)$ can have arbitrarily large dimension as an algebraic group: for example, if $S = \C[x_0,x_1,x_2]$ with weights $a$, $1$, and $1$, respectively, then $\dim(\mathrm{Aut}(S)) = \dim(S_1) + \dim(S_1) + \dim(S_a) = a + 6$.  Despite this, we'll prove that the Jordan constant of $\mathrm{Aut}(S)$ is uniformly bounded among polynomial rings $S$ with a fixed number of variables.  

Following the notation used in \Cref{diagonalizable}, we let $\mathcal{B} \coloneqq \{b : b = a_i \text{ for some } i\}$ be the set of positive integers that occur as a weight of the polynomial ring $S$.  For each $b \in \mathcal{B}$, $N_b$ is the number of weights equal to $b$.  Recall that $\bar{J}_N \coloneqq \bar{J}(\mathrm{GL}_N(\C))$.

\begin{lemma} \label{Jordan_const}
Let $S = \C[x_0,\ldots,x_{n+1}]$ be a weighted polynomial ring with weights $a_0,\ldots,a_{n+1}$. Then $\bar{J}(\mathrm{Aut}(S)) = \prod_{b \in \mathcal{B}} \bar{J}_{N_b}$.  In particular, for any integer $n$, there is a uniform upper bound $C_n$ on the weak Jordan constants of all groups $\mathrm{Aut}(S)$ where $S$ has $n+2$ variables.
\end{lemma}

\begin{proof}
Inside each graded piece $S_b$ with $b \in \mathcal{B}$, there is a subspace $V_b$ of dimension $N_b$ spanned by the variables of weight $b$, and a complementary subspace $W_b$ spanned by the remaining monomials of weighted degree $b$.  The direct sum $\bigoplus \mathrm{GL}(V_b)$ embeds as a subgroup of $\mathrm{Aut}(S)$, consisting of all automorphisms that don't ``mix" variables of different weights.  We'll show that any finite group $G \subset \mathrm{Aut}(S)$ is conjugate to a subgroup of $\bigoplus \mathrm{GL}(V_b) \subset \mathrm{Aut}(S)$. To do this, we'll construct the necessary coordinate change inside each $S_b$.

Since finite groups are linearly reductive in characteristic zero, the representation of $G$ on $S_b$ splits into a direct sum of irreducible representations.  In particular, since $W_b$ is $G$-invariant, we can find a complementary $G$-invariant subspace $V'_b$ inside $S_b$.  Define the change of coordinates on variables of weight $b$ in such a way that the span of the variables of weights $b$ becomes $V'_b$.  We can construct this change of coordinates independently within each $S_b$ and arrive at an automorphism of the entire graded ring $S$.  By construction, elements $g \in G$ don't mix variables of different weights in the new coordinates.

This proves that any finite group $G$ that appears as a subgroup of $\mathrm{Aut}(S)$ also embeds in $\bigoplus \mathrm{GL}(V_b)$.  Therefore,
$$\bar{J}(\mathrm{Aut}(S)) = \bar{J} \left( \bigoplus_{b \in \mathcal{B}} \mathrm{GL}(V_b) \right) = \prod_{b \in \mathcal{B}} \bar{J}_{N_b}.$$
Here the last equality comes from the general fact that $\bar{J}(G_1 \times G_2) = \bar{J}(G_1) \bar{J}(G_2)$ (this is one convenient property of weak Jordan constants that doesn't hold for regular Jordan constants).  Because there are $n+2$ weights total, we have $\sum_{b \in \mathcal{B}} N_b = n+2$. There are only finitely many possibilities for the collection of positive integers $N_b$ for a fixed $n$, so there is a uniform upper bound $C_n$ on the weak Jordan constant of $\mathrm{Aut}(S)$ depending only on $n$.  
\end{proof}

Since $\oplus_{b \in \mathcal{B}} \mathrm{GL}(V_b)$ is always a subgroup of $\mathrm{GL}_{n+2}(\C)$ via the block diagonal embedding, we may take $C_n = \bar{J}_{n+2}$ as the optimal upper bound in \Cref{Jordan_const}. 
 However, the exact values of the weak Jordan constants of $\mathrm{GL}_{N}(\C)$ are not known for most $N$, so the best known explicit value for $C_n$ in most cases is the (regular) Jordan constant $J_{n+2}$ computed by Collins \cite{Collins} (since $\bar{J}_N \leq J_N$).

\noindent {\bf Step 3}: Using \Cref{Jordan_const}, we now only need to bound the order of abelian subgroups of automorphisms $A \subset \mathrm{Aut}(S)$, where $S$ has $n+2$ variables.  Suppose that $A \subset H$ is an abelian subgroup of smallest index and assume we've changed coordinates on $S$ so that the action of $A$ is diagonal, using \Cref{diagonalizable}.

Now suppose that $f$ is a sum of $s$ monomials with nonzero coefficients and write it as
$$f = \sum_{i=1}^s K_{i} \prod_{j=0}^{n+1} x_j^{m_{ij}},$$
where each $K_i$ is nonzero by assumption. Package the exponents $m_{ij}$ into an $s \times (n+2)$ matrix $M$.  Each row corresponds to a monomial in $f$.  We can use \Cref{monomialexistence} to pick a distinguished collection of $n+2$ monomials in $f$: indeed, for each $j$, select a monomial of the form $x_j^{b_j}$ or $x_j^{b_j} x_k$, $k \neq j$ which has a nonzero coefficient in $f$.  Take only the $n+2$ rows of $M$ corresponding to these, and assemble them into a square $(n+2) \times (n+2)$ minor $B$ of $M$ in such a way that the monomial associated to $x_j$ goes in the $j$th row.  

\begin{lemma} \label{determinant}
The matrix $B$ constructed above is invertible and has determinant satisfying
$$0 < \mathrm{det}(B) \leq \frac{d^{n+2}}{a_0 \cdots a_{n+1}}.$$
\end{lemma}

\begin{proof}
We first note the following properties of $B$: first, every entry $b_j$ on the main diagonal is a positive integer satisfying $2 \leq b_j \leq d/a_j$. (The lower bound is by the criterion in \Cref{finite}, while the upper bound is because $x_j^{b_j}$ or $x_j^{b_j}x_k$ is a monomial of weighted degree $d$.) Second, each row of $B$ contains at most one nonzero element off the main diagonal; if it does, this element must be a $1$.

We can begin to compute the determinant by expanding along any rows or columns that have only one nonzero entry, namely the diagonal entry.  At each such step, the diagonal entry $b_j$ is a positive integer which is at most $d/a_j$, where $j$ is the index of the row in question.  After removing the $j$th row and column, the resulting minor always has the same properties as $B$, so it suffices to prove the inequality in the lemma with one copy of $d$ in the numerator and the $a_j$ in the denominator removed.  Continuing in this way, we may assume $B$ has exactly one off-diagonal $1$ in each row and column.  Up to a permutation of the indices, we can now further assume that $B$ is block diagonal with blocks of the form
$$\begin{pmatrix}
b_0 & 1 & 0 & \cdots & 0 \\
 0 & b_1 & 1 & \cdots & 0 \\
\vdots & \vdots & \ddots & \ddots & \vdots \\
0 & 0 & \cdots & b_{r-1} & 1 \\
1 & 0 & \cdots & 0 & b_r
\end{pmatrix}.$$
It now suffices to prove the lemma in the case that $B$ is a single block of the form above (so $r = n+1$). It's straightforward to compute that this ``loop matrix" has determinant $b_0 \cdots b_{n+1} + (-1)^{n+1} \neq 0$, so it is invertible (here we use $b_j \geq 2$).  As for the bound on the determinant, it automatically holds when $n$ is even since $b_0 \cdots b_{n+1} -1 < b_0 \cdots b_{n+1}$ and each $b_j < \frac{d}{a_j}$.  When $n$ is odd, use the series of equations $b_0 = (d-a_1)/a_0, b_1 = (d-a_2)/a_1, \ldots, b_{n+1} = (d-a_0)/a_{n+1}$ to compute
\begin{align*}
\det (B) & = b_{0} \cdots b_{n+1} + 1 =
\frac{(d-a_{0})(d-a_{1}) \cdots (d-a_{n+1})}{a_{0} \cdots a_{n+1}} + 1 \\
& = \frac{d(d-a_1)\cdots(d-a_{n+1}) - a_0(d-a_1)\cdots(d-a_{n+1}) + a_0 \cdots a_{n+1}}{a_0 \cdots a_{n+1}}
\\
& \leq
\frac{d(d-a_1)\cdots(d-a_{n+1})}{a_{0} \cdots a_{n+1}}
< \frac{d^{n+2}}{a_{0} \cdots a_{n+1}}.
\end{align*}
Here the penultimate inequality holds because $d-a_j \geq a_j$ for all $j$ by \Cref{finite}.
\end{proof}

With these properties of $B$ in hand, we return to the proof of the theorem.  We'll show that for any diagonal automorphism $x_j \mapsto c_j x_j$ in $A$, the scalars $c_j$ satisfy $|c_j| = 1$.  Indeed, since this automorphism preserves $f$, it preserves each monomial individually, and
$$K_{i} \prod_j (c_j x_{j})^{m_{ij}} = K_{i} \prod_j x_j^{m_{ij}},$$
for each $i = 1,\ldots,s$.  Therefore, $\prod_j c_j^{m_{ij}} = 1$.  Taking logarithms of the $|c_j|$, this means that $( \log |c_0|, \ldots, \log |c_{n+1}|) \in \ker M$.  But $\ker M \subset \ker B = \{0\}$ since $B$ is invertible, so $|c_j| = 1$ for each $j$.  Therefore, any element of $A$ can be represented as an $(n+2)$-tuple $(\theta_0,\ldots,\theta_{n+1})$ of elements of $\Q/\Z$, where $c_j = e^{2 \pi i \theta_j}$.

The condition that $(\theta_0,\ldots,\theta_{n+1})$ preserves $f$ can be expressed as $M (\theta_0,\ldots,\theta_{n+1})^\mathsf{T} \in \Z^s$.  We can obtain an upper bound for the order of $A$ by considering the weaker condition $B (\theta_0,\ldots,\theta_{n+1})^\mathsf{T} \in \Z^{n+2}$ instead (this says that at least the $n+2$ selected monomials in $f$ are preserved by the automorphism).  The number of distinct solutions for $(\theta_0,\ldots,\theta_{n+1})$ modulo $\Z^{n+2}$ to this latter equation is the index of $\Z^{n+2}$ in the superlattice spanned by $B^{-1} e_0, \ldots, B^{-1} e_{n+1}$, where $e_0,\ldots,e_{n+1}$ are the standard basis vectors in $\Z^{n+2}$.  This index equals $\det(B)$, so 
$$|A| \leq |\det(B)| \leq \frac{d^{n+2}}{a_0 \cdots a_{n+1}}$$ 
by \Cref{determinant}.  The original group $H$ which contained $A$ as a smallest index abelian subgroup therefore has order 
$$|H| \leq \bar{J}(\mathrm{Aut}(S))\frac{d^{n+2}}{a_0 \cdots a_{n+1}} \leq C_n \frac{d^{n+2}}{a_0 \cdots a_{n+1}}.$$
Finally, we have the desired bound 
$$|\on{Lin}(X)| = |G| = \frac{|H|}{d} \leq C_n \frac{d^{n+1}}{a_0 \cdots a_{n+1}}.$$
\end{proof}

\begin{example}[Fermat Hypersurfaces]
For any positive integer $n$ and degree $d \geq 3$, the {\it Fermat hypersurface} of dimension $n$ and degree $d$ in $\mathbb{P}^{n+1}$ is $X \coloneqq \{x_0^d + \cdots + x_{n+1}^d = 0\} \subset \mathbb{P}^{n+1}$.  Then $\on{Lin}(X)$ contains a copy of the symmetric group $S_{n+2}$ acting by permutation of the variables $x_i$ and the diagonal automorphisms given by multiplying each $x_i$ by an arbitrary $d$th root of unity (modulo the scalar transformations).  Therefore, $|\on{Lin}(X)| \geq (n+2)! d^{n+1}$. 

In fact, a computation of Shioda \cite{Shioda} shows that $|\on{Lin}(X)| = (n+2)! d^{n+1}$ when $X$ is defined over an algebraically closed field of characteristic zero (see also \cite{Kontogeorgis} for another proof and some generalizations of this result).  Note that a Fermat hypersurface can have extra automorphisms in positive characteristic \cite[Section 1]{Kontogeorgis}.  
\end{example}

This example shows that the order of growth with respect to the degree $d$ in the estimate \Cref{bound} is optimal and that we must have $C_n \geq (n+2)!$ for all $n$.  A natural question is whether we may actually take $C_n = (n+2)!$ for all $n$. We'll show that this is nearly true for $n = 1$, but not quite.

\begin{proposition} \label{curvebounds}
Let $X$ be a well-formed quasismooth weighted projective curve of degree $d$ in $\mathbb{P}(a,b,c)$ with finite linear automorphism group.  Then 
$$|\mathrm{Lin}(X)| \leq \frac{6d^2}{abc},$$
unless $a = b = c = 1$ and $X$ is projectively equivalent to one of the following two plane curves in $\mathbb{P}^2 = \mathbb{P}(1,1,1)$:
\begin{enumerate}
    \item The {\it Klein quartic} $xy^3 + yz^3 + zx^3 = 0$, with automorphism group isomorphic to $\mathrm{PSL}_2(\mathbb{F}_7)$ of order $168$;
    \item The {\it Wiman sextic} $10x^3 y^3 + 9x^5 z + 9y^5 z - 45x^2 y^2 z^2 - 135 xyz^4 + 27 z^6 = 0$, with automorphism group isomorphic to $A_6$ of order $360$.
\end{enumerate}
\end{proposition}

\begin{proof}
In order for a weighted projective curve to be well-formed, we must have that the weights $a, b$ and $c$ are pairwise relatively prime, and that each weight divides $d$.  Indeed, if some weight does not divide $d$, then the intersection $X \cap \mathbb{P}_{\mathrm{sing}}$ contains the corresponding coordinate point, which has codimension $1$ in $X$.  There are three cases to consider, depending on which of $a$, $b$, or $c$ coincide.

Suppose first that $a$, $b$, and $c$ are all distinct.  Then, any finite subgroup of $\mathrm{Aut}(\mathbb{P}(a,b,c))$ is abelian.  This follows from \Cref{Jordan_const}, which shows that $\bar{J}(\mathrm{Aut}(S)) = 1$ in this case.  Abelian subgroups of $\mathrm{Aut}(S)$ fixing the defining polynomial $f$ of $X$ have order at most $d^3/(abc)$ by {\bf Step 3} of the proof of \Cref{bound}; hence by {\bf Step 1}, $|\mathrm{Lin}(X)| \leq d^2/(abc)$.

Now suppose that $b = c$, but that $a$ is distinct from the other two weights.  Since $\mathbb{P}(a,b,c)$ is well-formed, we must have $b = c = 1$.  The weak Jordan constant of $\mathrm{Aut}(S)$ is $\bar{J}_1 \bar{J}_2  = \bar{J}_2 = 12$ in this case \cite[Section 2.2]{PS}.  Suppose our hypersurface is given by $X = \{f = 0\}$, where $f$ has weighted degree $d$.  In order for $|\mathrm{Lin}(X)|$ to exceed $6d^2/(abc) = 6d^2/a$, we would have to have (after conjugation) a finite subgroup $G$ of $\mathrm{GL}_1(\C) \oplus \mathrm{GL}_2(\C) \subset \mathrm{Aut}(S)$ fixing $f$ of order exceeding $6d^3/a$.  The maximal possible order of an abelian subgroup preserving $f$ is $d^3/a$, so we require our hypothetical $G$ to have no abelian subgroup of index less than or equal to $6$.  The image of $G$ under the projection
$$\mathrm{GL}_1(\C) \oplus \mathrm{GL}_2(\C) \xrightarrow{p_2} \mathrm{GL}_2(\C)$$
would also have no abelian subgroup of index less than or equal to $6$.  All finite subgroups of $\mathrm{GL}_2(\C)$ are central extensions of cyclic groups, dihedral groups, $A_4$, $S_4$, or $A_5$.  Of these, only $A_5$ has the required property: the largest abelian subgroup has index $12$.  Therefore, the image of $G$ in $\mathrm{PGL}_2(\C)$ is isomorphic to $A_5$.  It follows that $p_2(G)$ is a central extension of $A_5$ in $\mathrm{GL}_2(\C)$.  Since $X$ is quasismooth, the polynomial $f$ is of the form
$$f = x^{d/a} + x^{d/a-1}g_a(y,z) + x^{d/a-2}g_{2a}(y,z) +  \cdots + g_d(y,z),$$
for some polynomials $g_a, g_{2a}, \ldots, g_{d}$ of the indicated degrees in $y$ and $z$.  Here $g_d(y,z)$ must be nonzero since $f$ is irreducible. Each of the terms must be individually preserved by the action of $G$ because that action is block diagonal; in particular, $g_d(y,z)$ is an invariant polynomial under the action of $p_2(G)$.  But this means that the intersection of $p_2(G)$ with the center of $\mathrm{GL}_2(\C)$ has order at most $d$ (primitive roots of unity of higher degree could not preserve this polynomial), so $|p_2(G)| \leq |A_5|d = 60d$.  Similarly, $|\ker(p_2) \cap G| \leq d/a$, so
$$|G| \leq \frac{60d^2}{a}.$$
The combination of the inequalities $|G| \leq 60d^2/a$ and $|G| > 6d^3/a$ means that $d < 10$.  However, we've already seen that $g_d(y,z)$ is a polynomial invariant of the action of the binary icosahedral group in $\mathrm{GL}_2(\C)$.  The homogeneous generators for that invariant ring have degrees $12$, $20$, and $30$ by a result of Klein \cite{Klein2}, contradicting the bound on the degree.  It follows that no weighted projective curve of this form has more than $6d^2/a$ linear automorphisms.

The last possibility is that $a = b = c$, so that $X$ is a smooth plane curve. In this case, the problem of finding the largest possible automorphism groups in different degrees is well studied.  For $d \geq 4$, recall that all plane curve automorphisms are linear \cite{Chang}.  Klein \cite{Klein} computed the linear automorphism group of the quartic curve in \Cref{curvebounds}; it has the largest possible automorphism group of any curve of genus $g = 3$ by the Hurwitz bound $|\mathrm{Aut}(X)| < 84(g-1)$.  Wiman \cite{Wiman} first computed that the sextic in the proposition has automorphism group $A_6$.  Later work showed that the Wiman sextic is the unique degree $6$ curve with largest automorphism group up to projective equivalence \cite{DIK} and that the Fermat curve has the same property for various other $d \leq 20$ \cite{KMP1,KMP2}.  Finally, Harui \cite[Theorem 2.5]{Harui} proved that the two curves listed in \Cref{curvebounds} are the only ones with $|\mathrm{Lin}(X)| > 6d^2$ for {\it any} degree $d$.  This proves the proposition.
\end{proof}

This classification shows that we may take $C_1 = \frac{21}{2}$ in \Cref{bound}.  The author is unaware of any counterexamples to the theorem with $C_n = (n+2)!$ for $n \geq 2$.  By analogy with Collins' computations of Jordan constants, we might expect that unusual behavior such as in the $n = 1$ case occurs only for small $n$.

\begin{question} \label{Cnquestion}
Does \Cref{bound} hold with $C_n = (n+2)!$ for $n \geq 2$?  In particular, does the Fermat hypersurface have the largest automorphism group of any smooth hypersurface of degree $d$ in $\mathbb{P}^{n+1}$ when $n \geq 2$ and $d \geq 3$?
\end{question}

Many partial results in this direction are known for smooth hypersurfaces in ordinary projective space.  For instance, we have a fairly complete picture of the possible {\it orders} of automorphisms that can occur \cite{GAL,Zheng}.  The possible automorphism groups of smooth cubic surfaces over an algebraically closed field of arbitrary characteristic were classified by Dolgachev and Duncan \cite{DD}.  Moreover, the linear automorphism groups of smooth cubic threefolds and smooth quintic threefolds were classified by works of Wei and Yu \cite{WY} and Oguiso and Yu \cite{OY}, respectively.  The Fermat cubic fourfold is also known to have the largest possible automorphism group by a result of Laza and Zheng \cite{LZ}. In summary, the second part of \Cref{Cnquestion} is known to have an affirmative answer at least for the following pairs $(n,d)$: $(2,3), (3,3), (3,5),$ and $(4,3)$.

\section{Automorphisms of a Very General Hypersurface} \label{verygeneral}

Another result of Matsumura and Monsky is that the automorphism group of a {\it very general} hypersurface in $\mathbb{P}^{n+1}$ with $n \geq 2$ and degree $d \geq 3$ is trivial \cite[Theorem 5]{MM}.  One might hope that the same result holds for weighted projective hypersurfaces whenever the conditions of \Cref{finite} are met.  This turns out to be false, but we have the following result under a slightly stronger assumption on the degree.

\begin{theorem}
\label{generic}
Suppose that there exists a hypersurface $X \subset \mathbb{P}(a_0,\ldots,a_{n+1})$ of degree $d$ which is quasismooth and well-formed, where $n \geq 1$ and $d \geq 5 \max\{a_0,\ldots,a_{n+1} \}$.  Then for a very general such $X$, $\on{Lin}(X)$ is contained in the center of $\mathrm{Aut}(\mathbb{P})$ and is toric.  In particular, $\on{Lin}(X)$ is abelian.
\end{theorem}

\begin{proof}
Fix a very general hypersurface $X = \{f = 0\}$ with the given weights and degree. Any element of $\on{Lin}(X)$ comes from an automorphism $\alpha: S \rightarrow S$ of the graded ring $S = \C[x_0,\ldots,x_{n+1}]$.  The fact that $\alpha$ descends to $X$ means that $\alpha(f) = cf$ for some constant $c \in \C$.

The conditions of \Cref{generic} on the weights and degree are strictly stronger than those of \Cref{finite}, so we know that $\mathrm{Lin}(X)$ is finite.  In particular, the automorphism $\alpha$ has finite order.  It follows from \Cref{diagonalizable} that after conjugating by some automorphism of the graded ring $S$, $\alpha$ becomes diagonal, i.e., maps each $x_i$ to a scalar multiple of itself.  Let $\gamma: S \rightarrow S$ be such an automorphism that brings $\alpha$ into diagonal form.  Define $\beta \coloneqq \gamma \alpha \gamma^{-1}$ and $g \coloneqq \gamma (f)$, so that $\beta (g) = cg$.  For each $i$, let $c_i$ be the scalar such that $\beta(x_i) = c_i x_i$.  Next, let $G \coloneqq \mathrm{Aut}(S)$ and $H \coloneqq C_G(\beta)$ be the centralizer of the element $\beta$ in $G$.

We will show that unless $G = H$, that is, unless $\beta$ is actually contained in the center of $\mathrm{Aut}(S)$, the fact that $\{g = 0\}$ has automorphism $\beta$ forces more than $\dim(G) - \dim(H)$ monomials of degree $d$ in the polynomial $g$ to vanish.  This would contradict the assumption that $f$ was originally chosen to be very general, since the space of degree $d$ polynomials with $\beta$ as an automorphism would have codimension greater than $\dim(G/H)$.  The homogeneous space $G/H$, in turn, is isomorphic to the orbit of $\beta$ under conjugation.  This is the same idea used in Matsumura and Monsky's proof \cite[Theorem 5]{MM} of the analogous fact for hypersurfaces in $\mathbb{P}^n$. (In that paper, they considered both diagonalizable and unipotent automorphisms; since we are working over $\C$ instead of an arbitrary algebraically closed field, we only need to consider the former type.)

We've already seen that the dimension of $G = \mathrm{Aut}(S)$ is $\dim(G) = \dim(S_{a_0}) + \cdots + \dim(S_{a_{n+1}})$.  To compute the dimension of the centralizer $H$, it suffices to compute the dimension of its Lie algebra; we may do this by seeing which infinitesimal transformations commute with $\beta$.  Indeed, let $\sigma: x \mapsto x + \epsilon z$ be such a transformation, where $z = (z_0,\ldots,z_{n+1})$ is an $(n+2)$-tuple of homogeneous polynomials with $z_i \in S_{a_i}$.  We have that 
$$\sigma \beta (x_i) = \sigma (c_i x_i) = c_i \sigma(x_i) = c_i (x_i + \epsilon z_i),$$
while 
$$\beta \sigma (x_i) = \beta (x_i + \epsilon z_i) = c_i x_i + \epsilon \beta(z_i).$$
Comparing these two equations, we have that $\beta$ and $\sigma$ commute if and only if $\beta(z_i) = c_i z_i$; that is, if and only if each monomial in $z_i$ is multiplied by $c_i$ when applying $\beta$.  Therefore, $\dim(H)$ is equal to the cardinality of the set of ordered $2$-tuples $(i,y)$ such that $i \in \{0,\ldots,n+1\}$ and $y$ is a monomial of degree $a_i$ with $\beta(y) = c_i y$.  We may also describe the dimension of the entire group $G = \mathrm{Aut}(S)$ in a similar way: it is just the size of the set of {\it all} $2$-tuples $(i,y)$ with $y$ a monomial of degree $a_i$.  Therefore, $\dim(G/H) = |\Gamma|$, where $\Gamma$ is the set
$$\Gamma \coloneqq \{(i,y): \beta(y) \neq c_i y\}.$$
If $\Gamma$ is empty, then $G = H$, which is what we want.  Assuming that it is nonempty, we will now exhibit a vanishing monomial in $g$ for each $(i,y) \in \Gamma$, plus exactly one extra.  This would show that the number of vanishing monomials is greater than $\dim(G/H)$, as required. 

We'll begin by finding one vanishing monomial for each $(i,y) \in \Gamma$. Since $X$ is quasismooth, \Cref{monomialexistence} guarantees that we may choose a monomial of degree $d$ of one of the following two forms: $x_i^k$ or $x_i^k x_j$, $i \neq j$, for some positive integer $k$.  By default, we'll always choose the form $x_i^k$ when $a_i$ actually divides $d$.  By the assumption on degree, we have $k \geq 5$ in either case (if $k < 5$ and $ka_i + a_j = d$ for some $i, j$, then we must have $k = 4$, $a_i = a_j$ by the assumption on degree, so $5a_i = d$ and we can choose $x_i^5$ instead).  

If $a_i$ divides $d$, consider the pair of monomials $\{ x_i^{k-1} y, x_i^{k-2} y^2 \}$.  Since $(i,y) \in \Gamma$, $\beta(y) = c_y y$ where $c_y \neq c_i$.  But then, our two monomials are multiplied by $c_i^{k-1} c_y$ and $c_i^{k-2} c_y^2$, respectively, under $\beta$.  These constants cannot be equal or else we would have $c_y = c_i$.  If both monomials had nonzero coefficients in $g$, that would contradict the fact that $\beta(g) = cg$.  Therefore, at least one monomial of the pair vanishes in $g$.  The same reasoning works for the pair $\{ x_i^{k-1} y x_j,x_i^{k-2} y^2 x_j \}$ in the event that we began with $x_i^k x_j$ of degree $d$ instead.

This argument exhibits exactly $|\Gamma| = \dim(G/H)$ distinct vanishing monomials in the polynomial $g$. They are all distinct because any two pairs of monomials chosen above are disjoint.  This follows from the fact that we can recover the pair $(i,y)$ uniquely from either monomial of the pair.  This works as follows: given a monomial $x^I$ belonging to the pair we created from $(i,y) \in \Gamma$, find an index $i'$ such that: (1) the corresponding weight $a_{i'}$ is maximal among all variables appearing in $x^I$ with exponent at least $2$, and (2) the exponent of $x_{i'}$ is itself maximal among variables with indices satisfying the first condition.  Examining the forms of the monomials we chose above, one can show that since $k-2 > 2$, the index $i'$ identified by this procedure must be unique and equal to $i$. If we have two elements $(i,y_1)$ and $(i,y_2)$ in $\Gamma$ with $y_1 \neq y_2$, it's clear that the chosen pairs of monomials are disjoint. Thus, $y$ is also uniquely determined.

The final step is to find just one extra monomial in $g$ that vanishes.  To do this, we'll make a slight modification to the list of pairs above, without breaking the disjointness property of the previous paragraph.  Since $\Gamma$ is nonempty, we can fix a particular element $(i,y) \in \Gamma$.  Depending on the properties of $i$ and $y$, we find two vanishing monomials associated to $(i,y)$ rather than just one as follows:
\begin{itemize}
    \item If we chose $x_i^k$ with degree $d$ (here $k \geq 5$) and $y$ is not equal to some other variable $x_{i'}$, replace the pair $\{ x_i^{k-1} y,x_i^{k-2} y^2 \}$ with the two pairs $\{ x_i^k,x_i^{k-1} y \}$ and $\{ x_i^{k-2} y^2, x_i^{k-3} y^3 \}$.  We can do the same modification when we have $x_i^k x_j$ of degree $d$ (and $k \geq 5$, $y$ not linear) instead. None of the new monomials we've introduced can repeat among the ones we previously found for other elements in $\Gamma$.  We may now find two vanishing monomials for $(i,y)$ instead of one.
    \item If $x_i^k$ has degree $d$ with $k \geq 5$ and $y = x_{i'}$, then we still replace the pair $ \{ x_i^{k-1} x_{i'}, x_i^{k-2}x_{i'}^2 \}$ with the two pairs $\{ x_i^k,x_i^{k-1} x_{i'} \}$ and $\{ x_i^{k-2} x_{i'}^2, x_i^{k-3} x_{i'}^3 \}$.  However, the latter pair overlaps with the one we found for $(i',x_i) \in \Gamma$ in the special case that $k = 5$.  To remedy the issue in this one case, also replace the pair $\{ x_i x_{i'}^4, x_i^2 x_{i'}^3 \}$ associated to $(i',x_i)$ with $\{ x_{i'}^5,x_i x_{i'}^4 \}$.  As before, the process would be the same if we had started with $x_i^k x_j$ of degree $d$ in the beginning; no repeats are introduced.
\end{itemize}
By contradiction, we've now shown that $G = H$ so that $\beta$ is in the center of $\mathrm{Aut}(S)$.  This means that $\alpha = \gamma^{-1} \beta \gamma = \beta$, and more generally that $\alpha$ is diagonal in any choice of coordinates.  The induced automorphism of $\mathbb{P}(a_0,\ldots,a_{n+1})$ is therefore always toric, as claimed.
\end{proof}

For hypersurfaces satisfying the condition $d \geq 5 \max \{a_0,\ldots,a_{n+1}\}$ in \Cref{generic}, the stronger statement that $\mathrm{Lin}(X) = \{1\}$ for $X$ very general is not always true.

\begin{example}
Consider the family of hypersurfaces of degree $180$ in $\mathbb{P}^3(36,31,30,25)$.  The general $X$ in this family is quasismooth and well-formed.  Furthermore, the weights and degree satisfy the hypothesis of \Cref{generic}. However, since the only monomial of degree $180$ involving the variable $x_0$ of weight $36$ is $x_0^5$, any quasismooth $X$ has a non-trivial automorphism of order $5$ given by $x_0 \mapsto \zeta x_0$ for $\zeta$ a primitive fifth root of unity.  As predicted by the theorem, this automorphism is in the center of $\mathrm{Aut}(\mathbb{P})$.
\end{example}

\begin{commentA}
\begin{example}
    Here's another example where the degree is not in the range of \Cref{generic} but there is still a non-trivial automorphism in the center of $\mathrm{Aut}(\mathbb{P})$.  Consider the family of hypersurfaces of degree $48$ in $\mathbb{P}(16,13,12,9)$.  The general $X$ in this family is quasismooth and well-formed, but the only monomial of degree $48$ involving $x_0$ is $x_0^3$.  Therefore, $X$ always has a non-trivial automorphism of order $3$ given by $x_0 \mapsto \zeta x_0$ for $\zeta$ a primitive third root of unity.
\end{example}
\end{commentA}

In a similar way, one can construct examples where $d$ is arbitrarily large relative to the maximum of the weights, but non-trivial automorphisms still exist for any quasismooth $X$.  By having multiple ``isolated" weights, generic automorphism groups can be made to contain any abelian group.

Further, in the range where the conditions of \Cref{finite} apply but those of \Cref{generic} do not (i.e., when the degree satisfies $2 \max\{a_0,\ldots,a_{n+1}\} \leq d < 5 \max\{a_0,\ldots,a_{n+1}\}$), there are many examples of hypersurfaces with generic automorphisms outside the center of $\mathrm{Aut}(\mathbb{P})$.  These show that \Cref{generic} is close to optimal.

\begin{example}[Hyperelliptic curves, revisited]
    We saw above that hyperelliptic curves of genus $g$ naturally embed via the canonical map as $X_{2g+2} \subset \mathbb{P}^2(g+1,1,1)$.  Conversely, the general hypersurface of this degree is a hyperelliptic curve; up to a transformation of weighted projective space, a general equation becomes $x_0^2 + f(x_1,x_2) = 0$ and gives a double cover of $\mathbb{P}^1$.  The hyperelliptic involution is given by $x_0 \mapsto -x_0$, which is a nontrivial automorphism of $\mathbb{P}^2(g+1,1,1)$ that descends to $X$.  Further, this automorphism is not in the center of $\mathrm{Aut}(\mathbb{P})$.

    Similar reasoning works for other families of hypersurfaces $X_d \subset \mathbb{P}(a_0,\ldots,a_{n+1})$ in higher dimensions whenever $a_0 = d/2$; $X$ always has the involution given by the double cover, and this involution is often not in the center of $\mathrm{Aut}(\mathbb{P})$.
\end{example}

\begin{example}
It's well known that any cubic plane curve has non-trivial linear automorphisms.  Indeed, any smooth complex cubic plane curve $X_3 \subset \mathbb{P}^2(1,1,1) = \mathbb{P}^2$ is projectively equivalent to a curve in {\it Hesse normal form}
$$x^3 + y^3 + z^3 = 3C xyz,$$
where $C \in \C$, $C^3 \neq 1$.  This result dates back at least to the late 19th century \cite[v.3, p.22]{Weber}.  The curve $X$ defined by this equation has a linear automorphism group of order at least $18$ (in fact, the order equals $18$ except when $C$ takes one of a handful of special values \cite[Corollary 3.10]{BM}). For general $C$, this group is generated by permutations of $x,y,z$ and the transformation $(x:y:z) \mapsto (x:\zeta y: \zeta^2 z)$ for $\zeta$ a primitive third root of unity.  It acts transitively on the nine flex points of the curve, and is furthermore non-abelian. Since $\mathrm{Aut}(\mathbb{P}^2) = \mathrm{PGL}_3(\C)$ is centerless, all non-trivial automorphisms in this group must be outside the center. 

In the world of weighted projective hypersurfaces, we can bootstrap this example to any number of dimensions by looking at families such as $X_{15} \subset \mathbb{P}^4(5,5,5,3,3)$.  The variables corresponding to weights of $3$ and $5$ never mix, so we can again change coordinates so that the equation $f$ defining $X$ assumes Hesse normal form in the first three variables.  Then, the transformations above are still automorphisms of $f$, leaving the variables of weight $3$ unchanged.
\end{example}

\begin{example}
Examples of generic non-central automorphisms with $d = 4 \max \{a_0,\ldots,a_{n+1}\}$ also exist because of the fact that any quartic ``hypersurface" in $\mathbb{P}^1$, that is, a collection of four general points, has nontrival linear automorphism group.  More precisely, let $p_1,p_2,p_3,p_4 \in \mathbb{P}^1$ be four general points.  Any automorphism of $\mathbb{P}^1$ must preserve the cross-ratio of these four points, and the stabilizer of the cross-ratio under the permutation action of $S_4$ on these points is the Klein four-group $K = \{\mathrm{id}, (12)(34),(13)(24),(14)(23)\}$.  Further, since $\mathrm{PGL}_2(\C)$ is three-transitive, there is an automorphism mapping $p_1$, $p_2$ and $p_3$ according to any permutation $\sigma \in K$; by cross-ratio considerations, it also acts as $\sigma$ on the fourth.  It follows that the subgroup of $\mathrm{PGL}_2(\C)$ preserving this set of four points is isomorphic to $K \cong \Z/2 \oplus \Z/2$.

We can pick matrices in $\mathrm{GL}_2(\C)$ which descend to these transformations and preserve the quartic equation in two variables defining the given set of four points. This construction allows us to find many positive-dimensional examples with non-central automorphisms.  For instance, let $X_{20} \subset \mathbb{P}^3(5,5,4,4)$ be very general.  Then $X$ has nontrivial automorphisms defined by the same linear transformations as above in the first two variables and the identity on the variables of weight $4$.
\end{example}

\end{document}